\documentclass[12pt] {article}
\usepackage{amsmath,amsthm}
\usepackage{amssymb,latexsym}
\title{Hard Lefschetz theorem for valuations and related questions of integral geometry.}
\date{}
\author{ Semyon Alesker
\\  { \normalsize Department of Mathematics, Tel Aviv University, Ramat Aviv}
 \\  { \normalsize 69978 Tel Aviv,
Israel }
\\ {\normalsize e-mail: semyon@post.tau.ac.il}}
\newcommand{\RR}{\mbox{\rm $~\vrule height6.5pt width0.5pt
depth0.3pt\!\!$R}}

\newcommand{\CC}{\mbox{\rm $~\vrule height6.5pt width0.5pt
depth0.3pt\!\!$C}}

\def\eps{\varepsilon}

\def\lam{\lambda}

\def\str{\longrightarrow}

\def\qed { Q.E.D. }

\swapnumbers
\newtheorem{theorem}{Theorem}[section]
\newtheorem{corollary}[theorem]{Corollary}
\newtheorem{lemma}[theorem]{Lemma}
\newtheorem{proposition}[theorem]{Proposition}
\newtheorem{claim}[theorem]{Claim}
\theoremstyle{definition}

\newtheorem{definition}[theorem]{Definition}
\newtheorem{remark}[theorem]{Remark}
\theoremstyle{proposition-definition}
\newtheorem{proposition-definition}[theorem]{Proposition-Definition}

 \def\ck{{\cal K}}

\def\pt{\partial}
\begin{document}
\maketitle
\begin{abstract}
We continue studying the properties of the multiplicative
structure on valuations. We prove a new version of the hard
Lefschetz theorem for even translation invariant continuous
valuations, and discuss problems of integral geometry staying
behind these properties. Then we formulate a conjectural analogue
of this result for odd valuations.
\end{abstract}
 \setcounter{section}{-1}
\section{Introduction.}

   In this article we continue studying
properties of the multiplicative structure on valuations
introduced in \cite{alesker-mult}. In \cite{alesker-jdg} we have
proven a certain version of the hard Lefschetz theorem for
translation invariant even continuous valuations. Its statement is
recalled in Section 4 of this article (Theorem 4.1). The main
result of this article is Theorem 2.1 where we prove yet another
version of it for even valuations which is more closely related to
the multiplicative structure. As a consequence, we obtain a
version of it for valuations invariant under a compact group
acting transitively on the unit sphere and containing the operator
$-Id$ (Corollary 2.5). Then in Section 3 we state a conjecture
which is an analogue of the hard Lefschetz theorem for odd
valuations.

It turns out that behind the hard Lefschetz theorem for even
valuations stay results about the Radon and the cosine transforms
on the Grassmannians (see the proof of Theorem 2.1 in this
article, and also the proof of Theorem 1.1.1 in \cite{alesker-jdg}
which is also a version of the hard Lefschetz theorem). One passes
from even valuations to functions on Grassmannians using the
Klain-Schneider imbedding (see Section 1). The case of odd
valuations turns out to be related to integral geometry on partial
flags which seems to be not well understood. On the other hand,
various integral geometric transformations on such spaces can be
interpreted sometimes as intertwining operators (or their
compositions) for certain representations of $GL_n(\RR)$, and thus
can be reduced to a question of representation theory. This point
of view was partly used in \cite{alesker-bernstein} in the study
of the cosine transform on the Grassmannians, and more explicitly
in \cite{alesker-alpha} in the study of the generalized cosine
transform.  It would be of interest to understand these problems
for odd valuations.

Note also that the connection between the hard Lefschetz theorem
for valuations and integral geometry turns out to be useful in
both directions. Thus in \cite{alesker-jdg} it was applied to
obtain an explicit classification of unitarily invariant
translation invariant continuous valuations.

For the general background on convexity we refer to the book by
Schneider \cite{schneider-book}. For the classical theory of
valuations we refer to the surveys by McMullen and Schneider
\cite{mcmullen-schneider} and McMullen \cite{mcmullen-survey}.

\section{Background.} In this section we recall some definitions
and known results. Let $V$ be a real vector space of finite
dimension $n$. Let $\ck(V)$ denote the family of convex compact
subsets of $V$.

\begin{definition}
1) A function $\phi :\ck(V) \to \CC$ is called a valuation if for
any $K_1, \, K_2 \in \ck(V)$ such that their union is also convex
one has
$$\phi(K_1 \cup K_2)= \phi(K_1) +\phi(K_2) -\phi(K_1 \cap K_2).$$

2) A valuation $\phi$ is called continuous if it is continuous
with respect the Hausdorff metric on $\ck(V)$.

3) A valuation $\phi$ is called translation invariant if
 $\phi (K+x)=\phi(K)$ for every $x\in V$ and every $K\in \ck(V)$.

4) A valuation $\phi$ is called even if $\phi(-K)=\phi(K)$ for
every $K\in \ck(V)$.

5) A valuation $\phi$ is called homogeneous of degree $k$ (or
k-homogeneous) if for every $K\in \ck(V)$ and every scalar $\lam
\geq 0$, we have $\phi(\lam \cdot K)=\lam ^k \phi(K)$.
\end{definition}

We will denote by $Val (V)$ the space of translation invariant
continuous valuations on $V$. Equipped with the topology of
uniform convergence on compact subsets of $\ck (V)$ it becomes a
Fr\'echet space. We will also denote by $Val_k(V)$ the subspace of
$k$-homogeneous valuations from $Val (V)$. We will need the
following result due to P. McMullen \cite{mcmullen-euler}.
\begin{theorem}[\cite{mcmullen-euler}]
The space $Val(V)$ decomposes as follows
$$Val(V)=\bigoplus_{k=0}^{n} Val_k(V)$$
where $n=\dim V$.
\end{theorem}
In particular note that the degree of homogeneity is an integer
between 0 and $n=\dim V$. It is known that $Val_0(V)$ is
one-dimensional and it is spanned by the Euler characteristic
$\chi$, and $Val_n(V)$ is also one-dimensional and is spanned by a
Lebesgue measure \cite{hadwiger-book}. One has a further
decomposition with respect to parity:
$$Val_k(V)=Val_k^{ev}(V)\oplus Val_k^{odd}(V),$$
where $Val_k^{ev}(V)$ is the subspace of even $k$-homogeneous
valuations, and $Val_k^{odd}(V)$ is the subspace of odd
$k$-homogeneous valuations.

 Let us recall the imbedding of the space of
valuations into the spaces of functions on partial flags
essentially due to D. Klain \cite{klain1}, \cite{klain2} and R.
Schneider \cite{schneider-simple}. It will be used in Section 2 to
reduce the hard Lefschetz theorem for even valuations to integral
geometry of Grassmannians.

Let us denote by $Gr_i(V)$ the Grassmannian of real linear
$i$-dimensional subspaces in $V$. For a manifold $X$ we will
denote by $C(X)$ (resp. $C^\infty(X)$) the space of continuous
(resp. infinitely smooth) functions on $X$. Assume now that $V$ is
a Euclidean space. Let us describe the imbedding of
$Val_k^{ev}(V)$ into the space of continuous functions
$C(Gr_k(V))$ which we call Klain's imbedding. For any valuation
$\phi \in Val_k^{ev}(V)$ let us consider the function on $Gr_k(V)$
given by $L\mapsto \phi(D_L)$ where $D_L$ denotes the unit
Euclidean ball inside $L$. Thus we get a map $Val_k^{ev}(V)\to
C(Gr_k(V))$. The nontrivial fact due to D. Klain \cite{klain2}
(and heavily based on \cite{klain1}) is that this map is
injective.

\def\fkv{\tilde F(V)}
\def\cfk{C^-(\fkv)}
Now we would like to recall the Schneider imbedding of
$Val_k^{odd}(V)$ into the space of functions on a partial flag
manifold. Let us denote by $\fkv$ the manifold of pairs
$(\omega,M)$ where $ M\in Gr_{k+1}(V)$, and $\omega \in M$ is a
vector of unit length. Let us denote by $C^-(\fkv)$ the space of
continuous functions on $\fkv$ which change their sign when one
replaces $\omega$ by $-\omega$. Let us describe the imbedding of
$Val_k^{odd}(V)$ into $\cfk$ (following \cite{alesker-gafa}) which
we call Schneider's imbedding since its injectivity is an easy
consequence of a non-trivial result due to R. Schneider
\cite{schneider-simple} about characterization of odd translation
invariant continuous valuations. Fix a valuation $\phi\in
Val_k^{odd}(V)$. Fix any subspace $M\in Gr_{k+1}(V)$. Consider the
restriction of $\phi$ to $M$. By a result of P. McMullen
\cite{mcmullen-80} any $k$-homogeneous translation invariant
continuous valuation $\psi$ on $(k+1)$-dimensional space $M$ has
the following form. There exists a function $f\in C(S(M))$ (here
$S(M)$ denotes the unit sphere in $M$) such that for any subset
$K\in \ck(M)$
$$\psi(K)=\int_{S(M)}f(\omega) dS_k(K,\omega).$$
Moreover the function $f$ can be chosen to be orthogonal to any
linear functional (with respect to the Haar measure on the sphere
$S(M)$), and after this choice it is defined uniquely. We will
always make such a choice of $f$. If the valuation $\psi$ is odd
then the function $f$ is also odd. Thus applying this construction
to $\phi|_M$ for any $M\in Gr_{k+1}(V)$ we get a map
$Val_{k+1}^{odd}(V)\to \cfk$ defined by $\phi \mapsto f$. This map
turns out to be continuous and injective (see \cite{alesker-gafa}
Proposition 2.6, where the injectivity is heavily based on
\cite{schneider-simple}).

Let us recall the definition of the Radon  transform on the
Grassmannians. The orthogonal group acts transitively on
$Gr_i(V)$, and there exists a unique $O(n)$-invariant probability
measure (the Haar measure).

The Radon transform $R_{j,i}: Gr_i(V)\to Gr_j(V)$ for $j<i$ is
defined by $(R_{j,i}f)(H)= \int_{F \supset H}  f(F) \cdot dF$.
Similarly for $j>i$ it is defined by $(R_{j,i}f)(H)= \int_{F
\subset H}  f(F) \cdot dF$. In both cases the integration is with
respect to the invariant probability measure on all subspaces
containing (or contained in) the given one. The Radon transform on
real Grassmannians was studied in \cite{gelfand-graev-rosu},
\cite{grinberg}.

Recall now the definition of the cosine and sine of the angle
between two subspaces. Let $E\in Gr_{i}(V), \, F\in Gr_{j}(V)$.
Assume that $i\leq j$. Let us call {\itshape cosine of the angle}
between $E$ and $F$ the following number:
$$|cos(E,F)|:=\frac{vol_i(Pr_F(A))}{vol_i(A)},$$
where $A$ is any subset of $E$ of non-zero volume, $Pr_F$ denotes
the orthogonal projection onto $F$, and $vol_i$ is the
$i$-dimensional measure induced by the Euclidean metric. (Note
that this definition does not depend on the choice of a subset
$A\subset E$). In the case $i\geq j$ we define the cosine of the
angle between them as cosine of the angle between their orthogonal
complements:
$$|cos(E,F)| :=|cos(E^\perp ,F^\perp )|.$$
(It is easy to see that if $i=j$ both definitions are equivalent.)

Let us call {\itshape  sine of the angle} between $E$ and $F$ the
cosine between $E$ and the orthogonal complement of $F$:
$$|sin(E,F)|:= |cos(E,F^\perp )|.$$
 The following properties are well known (and rather trivial):
$$|cos(E,F)|=|cos(F,E)|=|cos(E^\perp ,F^\perp)|,$$
$$|sin(E,F)|=|sin(F,E)|=|sin(E^\perp ,F^\perp)|,$$
$$0\leq |cos(E,F)|, \, |sin(E,F)| \leq 1.$$

For any $1\leq i, \, j \leq n-1$ one defines the cosine transform
$$T_{j,i}:C(Gr_{i}(V)) \to C(Gr_{j}(V))$$ as follows:
$$(T_{j,i}f)(E):= \int_{Gr_{i}(V)} |cos(E,F)| f(F) dF,$$
where the integration is with respect to the Haar measure on the
Grassmannian. The cosine transform was studied in
\cite{matheron1}, \cite{matheron2}, \cite{goodey-howard1},
\cite{goodey-howard2}, \cite{goodey-howard-reeder},
\cite{alesker-bernstein}.

Let us state some facts on the multiplicative structure on
valuations. Let us briefly recall a construction of multiplication
from \cite{alesker-mult}, Section 1. A measure $\mu$ on a linear
space $V$ is called polynomial if it is absolutely continuous with
respect to a Lebesgue measure, and the density is a polynomial.
Let $\mu$ and $\nu$ be two polynomial measures on $V$  (the case
of Lebesgue measures would be sufficient for the purposes of this
article). Let $A,\, B\in \ck(V)$. Consider the valuations
$\phi(K)=\mu(K+A)$, $\psi(K)=\nu(K+B)$ where $+$ denotes the
Minkowski sum of convex sets. Let $\Delta:V\hookrightarrow V\times
V$ denote the diagonal imbedding. Then the product of valuations
is computed as follows:
$$(\phi \cdot \psi)(K)=(\mu \boxtimes \nu)(\Delta (K)+ (A\times
B))$$ where $\mu \boxtimes \nu$ denotes the usual product measure.
Product of linear combinations of measures of the above form is
defined by distributivity. The product turns out to be well
defined and it is determined uniquely by the above expressions
since the valuations of the above form are dense in the space of
polynomial valuations with respect to some natural topology.

We will discuss now only the translation invariant case though in
\cite{alesker-mult} the product was defined on a wider class of
polynomial valuations. Actually in this article we will need only
the following result.
\begin{proposition}
Let $p_1:V\to W_1$ and $p_2:V\to W_2$ be surjective linear maps.
Let $\mu_1$ and $\mu_2$ be Lebesgue measures on $W_1$ and $W_2$
respectively. Consider on $V$ the valuations
$\phi_i(K):=\mu_i(p_i(K)),\, i=1,2$.  Then
$$(\phi_1\cdot \phi_2)(K)=(\mu_1\boxtimes \mu_2) ((p_1\oplus
p_2)(K))$$ where $p_1\oplus p_2:V\to W_1\oplus W_2$ is given by
$(p_1\oplus p_2)(v)=(p_1(v),p_2(v))$ and $\mu_1\boxtimes \mu_2$ is
the usual product measure on $W_1\oplus W_2$.
\end{proposition}
{\bf Proof.} First recall the following well known formula (which
can be easily checked).  Let $I$ be a segment of unit length in
the Euclidean space $V$ orthogonal to a hyperplane $H$. Then one
has $$\frac{\pt}{\pt \lambda}\big |_0 vol(K+\lam I)=vol_{n-1}(Pr_H
K)$$ where $Pr_H$ denotes the orthogonal projection onto $H$. Let
now $I_1,\dots , I_k$ be pairwise orthogonal unit segments in $V$.
Let $L$ be the orthogonal complement to their span. Thus $\dim
L=n-k$. By the inductive application of the above formula one
obtains
$$\frac{\pt^k}{\pt \lam_1 \dots \pt \lam_k}\big |_0 vol_n
(K+\sum_{j=1}^k\lam_j I_j)=vol_{n-k} (Pr_L K).$$

Let us return to the situation of our proposition. We may assume
that $W_i$, $i=1,2$, are subspaces of $V$ and $p_i$ are orthogonal
projections onto them. We may also assume that the measures
$\mu_i$ coincide with the volume forms induced by the Euclidean
metric. Let us fix $I_1^{(i)},\dots, I_{m_i}^{(i)},\, i=1,2,$
pairwise orthogonal unit segments in the orthogonal complement to
$W_i$. Then
$$\phi_i(K)= \frac{\pt ^{m_i}}{\pt \lam_1\dots \pt \lam_{m_i}}\big |_0
vol_n(K+\sum _{j=1}^{m_i}\lam_j I_j^{(i)}).$$ Hence using the
construction of the product described above we get
$$(\phi_1\cdot \phi_2)(K)=$$
$$\frac{\pt ^{m_1}}{\pt \lam_1\dots \pt \lam_{m_1}} \frac{\pt
^{m_2}}{\pt \mu_1\dots \pt \mu_{m_2}}\big |_0 vol_{2n}
\left(\Delta(K)+ \sum_{j=1}^{m_1}\lam_j(I_j^{(1)}\times
0)+\sum_{l=1}^{m_2} \mu_l(0\times I_l^{(2)})\right)=$$
$$vol_{2n-m_1-m_2} (Pr_{W_1\oplus W_2}\left(\Delta
(K))\right)=vol_{2n-m_1-m_2}((p_1\oplus p_2)(K)).$$ \qed

\section{Hard Lefschetz theorem for even valuations.}
To formulate our  main result let us recall a definition from the
representation theory. Let $G$ be a Lie group. Let $\rho$ be a
continuous representation of $G$ in a Fr\'echet space $F$. A
vector $v\in F$ is called $G$-{\itshape smooth} if the map $G\str
F$ defined by $g\longmapsto g(v)$ is infinitely differentiable. It
is well known (and easy to prove) that smooth vectors form a
linear $G$-invariant subspace which is dense in $F$. We will
denote it by $F^{sm}$. It is well known (see e.g. \cite{wallach})
that $F^{sm}$ has a natural structure of a Fr\'echet space, and
the representation of $G$ in $F^{sm}$ is continuous with respect
to this topology. Moreover $(F^{sm})^{sm}=F^{sm}$. In our
situation the Fr\'echet space is $F=Val(V)$ with the topology of
uniform convergence on compact subsets of ${\cal K}(V)$, and
$G=GL(V)$. The action of $GL(V)$ on $Val(V)$ is the natural one,
namely for any $g\in GL(V),\, \phi\in Val(V)$ one has
$(g(\phi))(K)=\phi(g^{-1}K)$.
\begin{theorem}
Let $0\leq i<n/2$. Then the multiplication by $(V_1)^{n-2i}$
induces an isomorphism $Val_i^{ev}(V)^{sm}\tilde \to
Val_{n-i}^{ev}(V)^{sm}$. In particular for $p\leq n-2i$ the
multiplication by $(V_1)^p$ is an injection from
$Val_i^{ev}(V)\hookrightarrow Val_{i+p}^{ev}(V)$.
\end{theorem}
\begin{remark}
We would like to explain the use of the name "hard Lefschetz
theorem". The classical hard Lefschetz theorem is as follows (see
e.g. \cite{griffiths-harris}). Let $M$ be a compact K\"ahler
manifold of complex dimension $n$ with K\"ahler form $\omega$. Let
$[\omega]\in H^2(M,\RR)$ be the corresponding cohomology class.
Then for $0\leq i<n$ the multiplication by $[\omega]^{2(n-i)}$
induces an isomorphism $H^i(M,\RR)\tilde{\to} H^{2n-i}(M,\RR)$.
\end{remark}

First recall that in \cite{alesker-mult}, Theorem 2.6, we have
shown that $(V_1)^j$ is proportional to $V_j$ with a non-zero
constant of proportionality. Recall also the Cauchy-Kubota formula
(see e.g. \cite{schneider-weil}):
$$V_k(K)=c\int _{E\in Gr_k(V)}vol_k(Pr_E K)dE.$$
\begin{lemma}
Let $F\in Gr_i(V)$. Let $\phi(K)=vol_i(Pr_F(K))$. Then the image
of $V_k\cdot \phi$ in $C(Gr_{k+i})$ under the Klain imbedding is
given by the function
$$g(L)=c\int_{Gr_{k+i}(V)\ni R \supset  F}|cos(L,R)|dR$$
where $dR$ is the (unique) $O(i)\times O(n-i)$-invariant
probability measure on the Grassmannian of $(k+i)$-subspaces in
$V$ containing $F$, and $c$ is a non-zero normalizing constant.
\end{lemma}
{\bf Proof.} Using the Cauchy-Kubota formula and Proposition 1.3
we have
$$(V_k\cdot \phi)(K)=c\int_{E\in Gr_k(V)}vol_{k+i} ((Pr_F\oplus
Pr_E)(K))dE.$$ Let $K=D_L$ be the unit ball in a subspace $L\in
Gr_{k+i}(V)$. Then the image of $V_k\cdot \phi$ in $Gr_{k+i}(V)$
under the Klain imbedding is
$$g(L)=(V_k\cdot \phi)(D_L)=c \int_{E\in Gr_k(V)} vol_{k+i}
((Pr_E\oplus Pr_F)(D_L))dE.$$
\begin{claim}
$$vol_{k+i}((Pr_E\oplus Pr_F)(D_L))=\kappa |cos(L,(E+F))|\cdot
|sin(E,F)|$$ where $\kappa$ is a non-zero constant.
\end{claim}
Let us postpone the proof of this claim and let us finish the
proof of Lemma 2.2. We get
$$g(L)=c'\int _{E\in Gr_k(V)}|cos(L,(E+F))|\cdot |sin(E,F)|dE.$$
For a fixed subspace $F\in Gr_i(V)$ let us denote by $U_F$ the
open dense subset of $Gr_k(V)$ consisting of subspaces
intersecting $F$ trivially. Clearly the complement to $U_F$ has
smaller dimension. We have the map $T:U_F\to Gr_k(V/F)$ given by
$T(E):=(E+F)/F$. Clearly $T$ commutes with the action of the
stabilizer of $F$ in the orthogonal group $O(n)$ (which is
isomorphic to $O(i)\times O(n-i)$). Let $m:Gr_k(V/F)\to
Gr_{k+i}(V)$ be the map sending a subspace in $V/F$ to its
preimage in $V$ under the canonical projection $V\to V/F$. Then in
this notation we have
$$g(L)=c' \int _{E\in U_F} |cos(L,(m\circ T)(E)|\cdot
|sin(E,F)|dE.$$

Let us consider the following submanifold $M\subset
Gr_k(V/F)\times Gr_k(V)$ given by
$$M=\{(R,E)|\, E\subset m(R)\}.$$
Then $U_F$ is isomorphic to $M$ via the map $T_1:U_F\to M$ defined
by $T_1(E)=(T(E),E)$. Clearly $T_1$ commutes with the natural
action of the group $O(i)\times O(n-i)$ which is the stabilizer of
$F$ in $O(n)$. We have also the projection $t:M\to Gr_k(V/F)$
given by $t(R,E)=R$. Then
$$g(L)=c' \int_{R\in Gr_k(V/F)}dR |cos(L,m(R))|\left[ \int_{E\in
Gr_k(t^{-1}(R))} |sin(E,F)|d\mu_R(E)\right]$$ where $\mu_R$ is a
measure on  $Gr_k (t^{-1}(R))$. Note that the integral in square
brackets in the last expression is positive and does not depend on
$R$ since the map $T=t\circ T_1$ commutes with the action of
$O(i)\times O(n-i)$. Hence $g(L)=c''\int_{R\in Gr_k(V/F)}
|cos(L,m(R))|dR$. Lemma 2.2 is proved. \qed

{\bf Proof} of Claim 2.3. We may assume that $E\cap F=0$. Set for
brevity $p:=Pr_E\oplus Pr_F$. Then $p$ factorizes as $p=q\circ
Pr_{E+F}$ where $q:E+F\to E\oplus F$ is the restriction of
$Pr_E\oplus Pr_F$ to the subspace $E+F$. Then we have
$$vol_{k+i}((Pr_E\oplus Pr_F)(D_L))=vol_{k+i}(Pr_{E+F}(D_L))\cdot
\frac{vol_{k+i}(q(D_{E+F}))}{vol_{k+i}(D_{E+F})}= $$
$$vol_{k+i}(D_L)\cdot |cos(L,(E+F))|\cdot
\frac{vol_{k+i}(q(D_{E+F}))}{vol_{k+i}(D_{E+F})}.$$ Let us compute
the last term in the above expression. Note that the dual map
$q^*:E\oplus F\to E+F$ is given by $q^*((x,y))=x+y$. Clearly
$$\frac{vol(q^*(D_{E\oplus F}))}{vol D_{E\oplus F}}=|sin (E,F)|.$$
But the left hand side in the last expression is equal to
$\frac{vol(q(D_{E+ F}))}{vol(D_{E+ F})}.$ Thus Claim 2.3 is
proved. \qed

We will need one more lemma.
\begin{lemma}
Let $\phi\in Val_i^{ev}(V)$ be a valuation given by
$$\phi(K)=\int _{F\in Gr_i(V)}f(F)vol_i(Pr_F K)dF$$
with $f\in C(Gr_i(V))$. Then the image of $V_k\cdot \phi$ in
$C(Gr_{i+k}(V)$ under the Klain imbedding is given by
$$g(L)=c T_{k+i,k+i}\circ R_{k+i,i}(f)$$
where $c$ is a non-zero normalizing constant.
\end{lemma}
{\bf Proof.} By Lemma 2.2 one has
$$g(L)=c\int_{F\in Gr_i(V)}dF f(F)\int_{Gr_{k+i}\ni R\supset
F}|cos(L,R)|dR=$$
$$c\int _{R\in Gr_{k+i}(V)}dR |cos(L,R)|\int_{F\in Gr_i(R)}f(F)dF=
cT_{k+i,k+i}\circ R_{k+i,i}(f).$$ \qed

{\bf Proof} of Theorem 2.1. Let us consider the Klain imbedding
$Val_l^{ev}(V)^{sm}\hookrightarrow C^{\infty}(Gr_l(V))$. In
\cite{alesker-bernstein} it was shown that the image of
$Val_l^{ev}(V)^{sm} $ is a closed subspace which coincides with
the image of  the cosine transform $T_{l,l}:C^{\infty}(Gr_l(V))\to
C^{\infty}(Gr_l(V))$. By Lemma 2.4, $V_{n-2i}\cdot
T_{i,i}(f)=T_{n-i,n-i}\circ R_{n-i,i}(f)$. It was shown in
\cite{gelfand-graev-rosu} that $R_{n-i,i}: C^{\infty}(Gr_i(V))\to
C^{\infty}(Gr_{n-i}(V))$ is an isomorphism. Hence the image under
the Klain imbedding of $V_{n-2i}\cdot Val_i^{ev}(V)^{sm}$
coincides with the image of $T_{n-i,n-i}:C^{\infty}(Gr_{n-i})\to
C^{\infty}(Gr_{n-i}) $ which is equal to $Val_{n-i}^{ev}(V)^{sm}$.
Hence the multiplication by $V_1^{n-2i}: Val_i^{ev}(V)^{sm}\to
Val_{n-i}^{ev}(V)^{sm}$ is onto. Let us check that the kernel is
trivial. Indeed this operator commutes with the action of the
orthogonal group $O(n)$. The spaces $Val_i^{ev}(V)^{sm}$ and
$Val_{n-i}^{ev}(V)^{sm}$ are isomorphic as representations of
$O(n)$ (see \cite{alesker-jdg}, Theorem 1.2.2). Since every
irreducible representation of $O(n)$ enters with finite
multiplicity (in fact, at most 1) then any surjective map must be
injective. \qed

Now let us discuss an application to valuations invariant under a
group. Let $G$ be a compact subgroup of the orthogonal group
$O(n)$. Let us denote by $Val^G(V)$ the space of $G$-invariant
translation invariant continuous valuations. Assume that $G$ acts
transitively on the unit sphere. Then it was shown in
\cite{alesker-adv} that $Val^G(V)$ is finite dimensional. Also it
was shown in \cite{alesker-jdg} (see the proof of Corollary 1.1.3)
that $Val^G(V)\subset Val^{sm}(V).$ We have also McMullen's
decomposition with respect to the degree of homogeneity:
$$Val^G(V)=\oplus_{i=0}^n Val_i^G(V)$$
where $Val_i^G(V)$ denotes the subspace of $i$-homogeneous
$G$-invariant valuations.
 Then it was shown in \cite{alesker-mult}, Theorem 0.9 that $Val^G(V)$ is  a
finite dimensional graded algebra (with grading given by the
degree of homogeneity) satisfying the Poincar\'e duality (i.e. it
is a so called Frobenius algebra). A version of the hard Lefschetz
theorem was given in \cite{alesker-jdg}. Now we would like to
state another version of it.

\begin{corollary}
Let $G$ be a compact subgroup of the orthogonal group $O(n)$
acting transitively on the unit sphere. Assume that $-Id\in G$.
Let $0\leq i<n/2$ where $n=\dim V$. Then the multiplication by
$(V_1)^{n-2i}$ induces an isomorphism $Val_i^G(V)\tilde \to
Val_{n-i}^G(V)$. In particular for $p\leq n-2i$ the multiplication
by $(V_1)^p$ induces an injection $Val_i^G(V) \hookrightarrow
Val_{i+p}^G(V)$.
\end{corollary}
This result follows immediately from Theorem 2.1.

\section{The case of odd valuations.}
In this section we state the conjecture about an analogue of the
hard Lefschetz theorem for odd valuations and discuss it briefly.

{\bf Conjecture.} {\itshape Let $0\leq i<n/2$. Then the
multiplication by $(V_1)^{n-2i}$ is a map $Val_i^{odd}(V)^{sm} \to
Val_{n-i}^{odd}(V)^{sm}$ with trivial kernel and dense image. In
particular for $p\leq n-2i$ the multiplication by $(V_1)^p$ is an
injection from $Val_i^{odd}(V)^{sm}\hookrightarrow
Val_{i+p}^{odd}(V)^{sm}$. }

Let us show that it is enough to prove only either injectivity or
density of the image of the multiplication by $(V_1)^{n-2i}$. This
is a particular case of the following slightly more general
statement.
\begin{proposition}
Let $L:Val_i^{odd}(V)^{sm} \to Val_{n-i}^{odd}(V)^{sm}$ be a
continuous linear operator commuting with the action of the
orthogonal group $O(n)$. Then it has a dense image if and only if
its kernel is trivial.
\end{proposition}
{\bf Proof.} Note that the spaces $Val_i^{odd}(V)^{sm}$ and
$(Val_{n-i}^{odd}(V)^*)^{sm}\otimes Val_n(V)$ are isomorphic as
representation of the full linear group $GL(V)$ by the Poincar\'e
duality proved in \cite{alesker-mult}. Let us replace all spaces
by their Harish-Chandra modules. All irreducible representations
of $O(n)$ are selfdual. Hence the subspaces of $O(n)$-finite
vectors in $Val_i^{odd}(V)^{sm}$ and $Val_{n-i}^{odd}(V)^{sm}$ are
isomorphic as $O(n)$-modules. Moreover each irreducible
representation of $O(n)$ enters into these subspaces with finite
multiplicity (since by \cite{alesker-gafa} both spaces are
realized as subquotients in so called degenerate principal series
representations of $GL(V)$, namely representations induced from a
character of a parabolic subgroup). Then it is clear that the
operator $L$ between the spaces of $O(n)$-finite vectors commuting
with the action of $O(n)$ is surjective if and only if it is
injective. \qed

\section{Hard Lefschetz theorem for even valuations from
\cite{alesker-jdg}.} In this short section we remind another form
of the hard Lefschetz theorem for even valuations as it was proven
in \cite{alesker-jdg}. Note that it also was heavily based on the
properties of the Radon and cosine transforms on Grassmannians.

Let us fix on $V$ a Euclidean metric, and let $D$ denote the unit
Euclidean ball with respect to this metric.  Let us define on the
space of translation invariant continuous valuations an  operation
$\Lambda$ of mixing with the Euclidean ball $D$, namely
$$ (\Lambda \phi) (K):= \frac{d}{d\eps} \big|_{\eps =0} \phi(K +\eps D)$$
for any convex compact set $K$. Note that $\phi(K +\eps D)$ is a
polynomial in $\eps \geq 0$ by McMullen's theorem
\cite{mcmullen-euler}. It is easy to see that the operator
$\Lambda$ preserves parity and decreases the degree of homogeneity
by one. In particular we have
$$\Lambda: Val_k^{ev} \to Val^{ev}_{ k-1}(V).$$
The following result is Theorem 1.1.1 in \cite{alesker-jdg}.
\begin{theorem}
Let $n\geq k>n/2$. Then $\Lambda ^{2k-n} : (Val_k^{ev})^{sm} \to
(Val^{ev}_{n-k}(V))^{sm}$ is an isomorphism. In particular
$\Lambda ^i : Val_k^{ev} \to Val^{ev}_{k-i}(V)$ is injective for
$1\leq i \leq 2n-k$.
\end{theorem}

\end{document}